\documentclass[12pt]{article}

\usepackage{amsmath}
\usepackage{amsthm,amssymb}
\usepackage{multirow}
\usepackage{slashbox}
\usepackage{mathrsfs}
\usepackage{float}
\usepackage{color}
\usepackage{graphicx}
\usepackage{subfigure}
\usepackage{tabu}
\usepackage{booktabs}
\usepackage{array}
\usepackage{caption}
\usepackage{multirow}
\usepackage{geometry}
\usepackage{verbatim}

\makeatletter
\@addtoreset{equation}{section}
\makeatother

\newtheorem{de}{Definition}[section]%定义
\newtheorem{thm}{Theorem}[section]%定理
\newtheorem{lem}{Lemma}[section]%引理
\newtheorem*{pr}{Proof}%无标号证明
\newtheorem{remark}{Remark}[section]

\newcommand{\A}{{\mathcal A}}
\newcommand{\B}{{\mathcal B}}
\newcommand{\C}{{\mathcal C}}
\newcommand{\D}{{\mathcal D}}
\newcommand{\I}{{\mathcal I}}
\newcommand{\Q}{{\mathcal Q}}
\newcommand{\U}{{\mathcal U}}
\newcommand{\s}{{\mathcal S}}
\newcommand{\V}{{\mathcal V}}

\newcommand{\M}{{\mathcal M}}

\newcommand{\Aj}{{\hat {\mathcal A}}}
\newcommand{\Bj}{{\hat {\mathcal B}}}

\newcommand{\nnl}{{n \times n \times \ell}}

\newcommand{\pql}{{p \times q \times \ell}}

\newcommand{\nn}{{n \times n}}

\newcommand{\RR}{{\mathbb R}}
\newcommand{\CC}{{\mathbb C}}

\newcommand{\fold}{{\mathtt {fold}}}
\newcommand{\unfold}{{\mathtt {unfold}}}
\newcommand{\bcirc}{{\mathtt {bcirc}}}
\newcommand{\diag}{{\mathtt {diag}}}
\newcommand{\fft}{{\mathtt {fft}}}
\newcommand{\ifft}{{\mathtt {ifft}}}

\newcommand{\T}{{\rm {T}}}

\newcommand{\Ajj}{{\hat {A}}}
\newcommand{\Bjj}{{\hat {B}}}

\newcommand{\sjj}{{\hat {S}}}

\title{A note on perturbation analysis for T-product based tensor singular values}
\author{
Yating Zhang\footnote{School of Mathematical Sciences, Ocean University of China, Qingdao 266100, China.
E-Mail: {\tt zhangyating@stu.ouc.edu.cn}},
Xiaoxia Guo\footnote{School of Mathematical Sciences, Ocean University of China, Qingdao 266100, China. E-Mail: {\tt guoxiaoxia@ouc.edu.cn}},
Pengpeng Xie\footnote{Corresponding author: School of Mathematical Sciences, Ocean University of China, Qingdao 266100, China.
E-Mail: {\tt xie@ouc.edu.cn}. The work of this author is supported in part by the NSFC grant 11801534.}.
}
\date{}          %%不显示日期

\begin{document}

\maketitle                %%显示题目
%%摘要、关键词
\vspace{-1.2 cm}
\begin{abstract}
In this note, we present perturbation analysis for the T-product based tensor singular values defined by Lu et al. First, the Cauchy's interlacing-type theorem for tensor singular values is given. Then, the inequalities about the difference between the singular values of two matrices proposed by Mirsky are extended to tensor cases. Finally, we introduce some useful inequalities for the singular values of tensor products and sums.
\\ \hspace*{\fill} 
\end{abstract}

%%===================写背景，发展史，应用等：Introduction==================
\section{Introduction}

%%%==================第一自然段：T-product + T-SVD引入===============
\hskip 2em In 2008, Kilmer, Martin and Perrone \cite{Kilmer08} proposed the tensor operation T-product and extended the matrix SVD to third-order tensors. Based on this tensor-tensor multiplication, many scholars have contributed to the theoretical results and algebraic analysis \cite{Kilmer13,MiaoQi20}. Specific applications of the T-product and the T-SVD can be seen in \cite{Semerci14,XiaoChen21,Yang16,Zhang17,Zhou18}.

%%=================第二自然段：介绍奇异值的定义：两个论文中分别提了========
%先提不同的定义，以及本质。然后介绍新定义的作用。再写旧定义的作用。最后我们用啥。
\hskip 2em Under the T-SVD framework, two different definitions of tensor singular values have been proposed successively by Lu et al. \cite{Lu20} and Qi and Yu \cite{Qi21}. 
The former definition can be understood as the average of the singular values of several matrices, 
while the latter definition which is called the T-singular value can be interpreted as a measure of the length of some singular value vectors. 
In \cite{Qi21}, the T-singular values are used to define the tail energy of a third-order tensor, 
and applied to the error estimation of a tensor sketching algorithm for low-rank tensor approximation. 
In \cite{Ling21}, Ling et al. presented four necessary conditions and a set of sufficient and necessary conditions for s-diagonal tensors via the T-singular values.
Contrastly, the main contribution of tensor singular values defined in \cite{Lu20} is to help define the tensor nuclear norm. 
The definition and properties of the nuclear norm are consistent with the matrix cases. 
Furthermore, the tensor nuclear norm can be used to solve tensor robust principal component analysis problem and applied to image recovery and background modeling problems. 
We hope to continue along the notion of \cite{Lu20} to discover more properties and effects of the tensor singular values.

%%%=================第三自然段：论文第三部分：奇异值交错定理研究动力============
%主要是特征值的定义两篇论文+近期那篇论文的研究贡献：关于特征值的有关不等式。
\hskip 2em Braman \cite{Braman10} introduced the eigendecomposition of $n \times n \times n$ tensors firstly. A recent work by Liu and Jin \cite{Liu21} studied the T-eigenvalues of third-order tensors in detail. They proved some T-eigenvalue inequalities for Hermitian tensors, including extensions of Weyl's theorem and Cauchy's interlacing theorem from the matrix case to the tensor case. Moreover, they described the stability of the T-eigenvalues and studied the Lyapunov equation for tensors. It is therefore natural to consider the extension of the classical perturbation theory for matrix singular value problems. We show herein that the interlacing property and the Mirsky-type inequalities can be recast in the T-product formalism. Also, some singular value inequalities for the tensor products and sums can be generalized.

%%%==========第四自然段：本文的研究内容+论文章节结构===================
\hskip 2em This paper is organized as follows. In section 2, we review basic definitions and notations. 
In section 3, we consider the singular values of subtensors.
 Generalized Mirsky inequalities comparing the differences of two tensors with the singular values of their difference are proposed in section 4. 
 Section 5 shows other useful inequalities for the singular values of the tensor products and sums. A conclusion is given in section 6.

%%============论文的第二部分：已经存在的定义定理等内容================
%%=============Preliminaries===========================================
\section{Preliminaries}

\subsection{Notations}
\hskip 2em The term tensor refers to a multidimensional array of numbers and it is necessary to break a tensor up into various slices and tubal elements, and to have an indexing on those. 
In our paper, tensors are denoted by calligraphic letters and matrices are denoted by capital letters.
% A third-order tensor $\A\in\CC^{m \times n \times p}$,
 We denote $\A^{(i)}\equiv\A(i,:,:)$, $\overrightarrow{\A_i}\equiv\A(:,i,:)$ and  $A_i\equiv\A(:,:,i)$ for the $i$-th 
 horizontal, lateral and frontal slice of a third-order tensor $\A\in\CC^{m \times n \times p}$ respectively. 
 We also use $a_{ijk}$ to represent its $(i,j,k)$-th element, and $\A_{i,j}\equiv\A(i,j,:)$ for its $\left(i,j\right)$-th tubal scalar.

\hskip 2em For a third-order tensor $\A$, as in \cite{Kilmer13}, it is defined that
\begin{equation*}
\unfold(\mathcal{A})=
\begin{bmatrix}
 A_1 \\
 \vdots        \\
 A_{p}
\end{bmatrix}\in\CC^{mp \times n},\quad
\fold(\unfold(\mathcal{A}))=\mathcal{A}.
\label{fold}
\end{equation*}

\hskip 2em The discrete Fourier transform (DFT) can transform block circulant matrices into block diagonal matrices.  Mathematically, this means that if $F_n=\omega^{(j-1)(k-1)}\in\CC^{n \times n}$ with $\omega=e^{-\frac{2\pi \mathtt{i}}{n}}$ denotes the DFT matrix, then we can obtain
\begin{equation*}
\left(F_{p} \otimes I_{m}\right)
\bcirc \left(\A\right)
\left(F_{p}^{-1} \otimes I_{n}\right)=\diag\left(\Ajj_1,\Ajj_2,\dotsc,\Ajj_{p}\right)
,
\label{block diagonal}
\end{equation*}
where $\otimes$ denotes the Kronecker product and $\bcirc(\mathcal{A})$ is block-circulant matrix of the form
\begin{equation*}
\bcirc(\mathcal{A})=
\begin{bmatrix}
A_1      &A_{p}    &A_{p-1}  &\dotsc       &A_2      \\
A_2      &A_1      &A_{p}    &\dotsc       &A_3      \\
\vdots   &\ddots   &\ddots   &\ddots       &\vdots   \\
A_{p}    &A_{p-1}  &\ddots   &A_2          &A_1
\end{bmatrix}.
\label{bcirc}
\end{equation*}

\hskip 2em Using Matlab  notations, $\Aj\equiv\fft\left(\A,[\ ],3\right)\in\CC^{m \times n \times p}$ with $\Aj(:,:,i)=\Ajj_i$ is the tensor obtained by applying the fast Fourier transform (FFT) along each tubal scalar of $\A$.
\begin{comment}Naturally, $\A\equiv\ifft(\Aj,[\ ],3)$ and it is not hard to show that
\begin{equation}
\Ajj_i    =\sum^{p}_{j=1}\omega^{(i-1)(j-1)} A_{j},\ i=1,2,\dotsc,p.
\label{relation}
\end{equation}
\end{comment}

\subsection{Definitions and propositions}
\hskip 2em In this subsection, we give the basic definitions and propositions from \cite{Kilmer08,Kilmer13,Lu20}.

%--------------------定义2.1：C上张量乘积------------------------
\begin{de} (T-product) Suppose $\A\in\CC^{m \times n \times p}$ and $\B\in\CC^{n \times t \times p}$, then the T-product $\A*\B$ is the tensor in $\CC^{m \times t \times p}$
\begin{equation*}
 \A*\B=\fold\left(\bcirc\left(\A\right)\cdot\unfold\left(\B\right)\right).
 \label{Tproduct}
\end{equation*}
\end{de}

%----------------------引理2.1:-----------------
\begin{lem} Suppose tensors $\A,\B$ and $\C$ are well-defined, then
\begin{equation*}
\C=\A\pm\B\iff\hat{C}_{i}=\hat{A}_i\pm\hat{B}_i,\quad \C=\A*\B\iff\hat{C}_{i}=\hat{A}_i\hat{B}_i.
\end{equation*}
\begin{comment}
\begin{equation*}\C=\A\pm\B\iff\hat{C}_{i}=\hat{A}_i\pm\hat{B}_i,\label{sum}\end{equation*}% 公式标号2.3
\begin{equation*}\C=\A*\B\iff\hat{C}_{i}=\hat{A}_i\hat{B}_i.\label{product}\end{equation*}% 公式标号2.4
\end{comment}
\end{lem}

\begin{comment}%块张量的定义直接在性质里介绍就行
%--------------------定义2.2：C上的块张量定义------------------
\begin{de} (Block Tensor) Suppose $\A\in\CC^{m_1 \times n_1 \times n}$, $\B\in\CC^{m_1 \times n_2 \times n}$, $\C\in\CC^{m_2 \times n_1 \times n}$ and $\D\in\CC^{m_2 \times n_2 \times n}$. The block tensor
\begin{equation*}
\begin{bmatrix}
 \A & \B \\
 \C & \D
\end{bmatrix}\in\CC^{(m_1+m_2) \times (n_1+n_2) \times n}
\end{equation*}
is defined by compositing the frontal slices of four tensors.
\end{de}
\end{comment}

\begin{comment}
%%------------------引理2.2：fft的先后问题，放在块张量定义后面--------------------
\begin{lem} Let $\A\in\RR^{m_1 \times n_1 \times n}$, $\B\in\RR^{m_1 \times n_2 \times n}$, $\C\in\RR^{m_2 \times n_1 \times n}$ and $\D\in\RR^{m_2 \times n_2 \times n}$, then
\begin{equation*}
 \rm{fft}\left(
 \begin{bmatrix}
 \A & \B \\
 \C & \D
 \end{bmatrix}
 ,[\ ],3 \right)=
 \begin{bmatrix}
  \rm{fft}\left(\A,[\ ],3\right) &
  \rm{fft}\left(\B,[\ ],3\right)\\
  \rm{fft}\left(\C,[\ ],3\right) &
  \rm{fft}\left(\D,[\ ],3\right)
 \end{bmatrix}
 \in\CC^{(m_1+m_2) \times (n_1+n_2) \times n},
\end{equation*}
where the tensor in the left is the block tensor defined by compositing the frontal slices of each tensors.
\end{lem}
\end{comment}

%-----------------------定义2.2:张量R上 转置---------------
\begin{de} (Tensor Transpose) Suppose $\A\in\RR^{m \times n \times p}$, then $\A^{\T}$ is the $n\times m\times p$ tensor obtained by transposing each of the frontal slices and then reversing the order of transposed slices 2 through $p$.
\end{de}

%-----------------------------定义2.3：R上单位张量--------------------------
\begin{de} (Identity Tensor) The $\nnl$ identity tensor $\I_{nn\ell}$ is the tensor whose first frontal slice is the $\nn$ identity matrix, and whose other frontal slices are all zeros.
\end{de}

%--------------------------定义2.4:实数域，逆张量-----------------
\begin{de} (Tensor Inverse) The tensor $\A\in\RR^{\nnl}$ has an inverse $\B$ provided that
\begin{equation*}
\A*\B=\B*\A=\I.
\end{equation*}
\end{de}

%-------------------------定义2.5:R上正交与部分正交张量----------------
\begin{de} (Orthogonal Tensor) The tensor $\Q\in\RR^{\nnl}$ is orthogonal if \begin{equation*}\Q^{\T}*\Q=\Q*\Q^{\T}=\I.\end{equation*}
The tensor $\Q\in\RR^{\pql}$ is partially orthogonal if \begin{equation*}\Q^{\T}*\Q=\I_{qq\ell}.\end{equation*}
\end{de}

%======================插入张量的F范数定义:定义2.6============================
\begin{de} (Tensor Frobenius Norm) Suppose $\A=\left(a_{ijk}\right)\in\RR^{m \times n \times p}$, then
\begin{equation*}
\parallel\A\parallel_F=\sqrt{\sum^{m}_{i=1}\sum^{n}_{j=1}\sum^{p}_{k=1} a_{ijk}^2}.
\end{equation*}
\end{de}

%=========================定理2.1：张量的F范数酉不变性========================
\begin{thm} If $\Q$ is an orthogonal tensor, then
\begin{equation*}
\parallel\Q*\A\parallel_F=\parallel\A\parallel_F.
\end{equation*}
\end{thm}

%%=======================定理2.2：FFT前后的两个张量的F范数关系等式============
\begin{thm} Suppose $\A\in\RR^{m \times n \times p}$ and $\Aj=\fft(\A,[\ ],3)$, then
\begin{equation}
\parallel\A\parallel_F=\frac{1}{\sqrt{p}}\parallel\Aj\parallel_F.
\label{Fnorm}
\end{equation}
\end{thm}

%\begin{comment}
%%----------------------------定义2.7:f对角张量定义------------------
%\begin{de} (F-diagonal Tensor) We say a tensor is f-diagonal if its each frontal slice is diagonal.
%\end{de}
%\end{comment}

%=====================================T-SVD定理=======================================
\begin{thm} (T-SVD) Let $\A\in\RR^{m\times n\times p}$. Then it can be factorized as
\begin{equation*}
\A=\U*\s*\V^{\T}
\end{equation*}
where $\U\in\RR^{m\times m\times p}$, $\V\in\RR^{n\times n\times p}$ are orthogonal, and $\s\in\RR^{m\times n\times p}$ is an f-diagonal tensor, $i.e.$, each frontal slice is diagonal.
\end{thm}

%=========================奇异值定理与相关性质说明======================================
\hskip 2em The entries on the diagonal of the first frontal slice $\s(:,:,1)$ of $\s$ have the decreasing  order property, $i.e.$,
\begin{equation*}
\s(1,1,1)\geq\s(2,2,1)\geq\cdots\geq\s(n',n',1)\geq0.
\end{equation*}
where $n'= \min(m,n)$. The above property holds since the inverse FFT gives
\begin{equation}
\s(i,i,1)=\frac{1}{p}\sum^{p}_{j=1}\hat{\s}(i,i,j).
\label{singualr}
\end{equation}
\hskip 2em The entries on the diagonal of $\hat{\s}(:,:,j)$ are the singular values of $\hat{\A}(:,:,j)$.

\begin{de} (Tensor singular values \cite{Lu20}) The entries on the diagonal of $\s(:,:,1)$ are the singular values of $\A$.
\end{de}

%%=====================第三部分：只写奇异值交错定理===================
\section{The singular values of subtensors}

\hskip 2em Liu and Jin \cite[Theorem 4.5]{Liu21} have extended Cauchy's interlacing
theorem from the matrix case to the tensor case. This theorem gives the inequality about the T-eigenvalues of one tensor and its subtensors. Similarly, in this section, we study the relationship between singular values of one tensor and its subtensors. In other words, our theorem is an extension of the matrix singular value interlacing theorem.

%%=========奇异值交错定理：直接写张量形式，========矩阵形式的在证明中体现=====
\begin{thm} Suppose the tensor $\A\in\RR^{m \times n \times r}$ with singular values
\begin{equation*}
\alpha_1\geq\alpha_2\geq\cdots\geq\alpha_{n'}\geq 0,\ n'=\min\{m,n\}.
\end{equation*}
Suppose $\B\in\RR^{p \times q \times r}$ is a subtensor of $\A$ with singular values
\begin{equation*}
\beta_1\geq\beta_2\geq\cdots\geq\beta_{p'}\geq 0,\ p'=\min\{p,q\}.
\end{equation*}
Then,
\begin{align}
\alpha_j&\geq\beta_j,\ j=1,2,\dotsc,p',\label{jiaocuo_result1}\\
\beta_j&\geq\alpha_{j+(m-p)+(n-q)},\ j=1,2,\dotsc,\min\{(p+q-m),(p+q-n)\}.\label{jiaocuo_result2}
\end{align}
\end{thm}

\begin{pr} \upshape{Suppose} every frontal slice $\Ajj_i\in\CC^{m \times n}$ of $\Aj=\fft(\A,[\ ],3)$ and $\Bjj_i\in\CC^{p \times q}$ of $\Bj=\fft(\B,[\ ],3)$ have the singular values numbered in decreasing order
\begin{equation*}
\sigma^{i}_1\geq\cdots\geq\sigma^{i}_{n'}\geq 0,\quad
\tau^{i}_1\geq\cdots\geq\tau^{i}_{p'}\geq 0,\quad
i=1,2,\dotsc,r.
\end{equation*}
From (\ref{singualr}), the singular values of the tensor $\A$ and subtensor $\B$ are
\begin{equation*}
\alpha_j=\frac{1}{r}\sum_{i=1}^r\sigma_j^i,\ j=1,2,\dotsc,n',
\beta_j=\frac{1}{r}\sum_{i=1}^r\tau_j^i,\ j=1,2,\dotsc,p'.
\end{equation*}
Since the singular value  interlacing theorem \cite[Theorem 1]{Thompson71} tells us that for any $i=1,2\dotsc,r$,
\begin{align}
\sigma^i_j&\geq\tau^i_j,\ j=1,2,\dotsc,p',\label{matrixresult1}\\
\tau^i_j&\geq\sigma^i_{j+(m-p)+(n-q)},\ j=1,2,\dotsc,\min\{(p+q-m),(p+q-n)\},\label{matrixresult2}
\end{align}
the results (\ref{jiaocuo_result1}) and (\ref{jiaocuo_result2}) follows from (\ref{matrixresult1}) and (\ref{matrixresult2}) respectively.\qedsymbol
\end{pr}

\begin{remark} \upshape{This theorem also suggests that if a lateral slice is added to a tensor $\A\in\RR^{m\times n\times r}$ with $m>n$, then the largest singular value increases and the smallest one decreases.} That is, if $\B=[\A,\overrightarrow{\M}]\in\RR^{m\times (n+1)\times r}$, then
\begin{equation*}
\sigma_{\max}(\B)\geq\sigma_{\max}(\A),\quad\sigma_{\min}(\B)\leq\sigma_{\min}(\A).
\end{equation*}
When $r=1$, this conclusion is consistent with standard matrix algebra operations and terminology \cite[Corollary 2.4.5]{Golub13}.
\end{remark}

%%%==========第四部分：只写Mirsky类型的两个不等式=================
%%%=============两个不等式写在一个定理里边=======================
\section{Mirsky-type singular value inequalities for tensors}
\hskip 2em For any unitarily invariant norm of matrices, it is proved by Mirsky \cite{Mirsky60} that
\begin{equation}
\parallel\diag(\sigma_1-\widetilde{\sigma}_1,\dotsc,\sigma_n-\widetilde{\sigma}_n)\parallel
\leq\parallel B-\widetilde{B}\parallel,
\label{matrixnorm}
\end{equation}
where $\sigma_i's$ and $\widetilde{\sigma}_i's$ are the corresponding singular values of $B$ and $\widetilde{B}$.
%which is proved by Mirsky \cite{Mirsky60}.
To the best of our knowledge, there is no straightford definition of unitarily invariant norm for tensors. 
Instead, we turn to the Frobenius norm and the spectral norm (also known as the largest singular value \cite{Lu20}) based tensor Mirsky-type inequalities.

\begin{thm} Suppose $\A,\B\in\RR^{m\times n\times r}$ with $n'=\min\{m,n\}$ and
\begin{equation*}
\alpha_1\geq\cdots\geq\alpha_{n'}\geq 0,\ \beta_1\geq\cdots\geq\beta_{n'}\geq 0
\end{equation*}
are the singular values of $\A$ and $\B$ respectively. Then,
\begin{equation}
\left[\sum_{i=1}^{n'}\left(\alpha_i-\beta_i\right)^2\right]^{\frac{1}{2}}\leq\parallel\B-\A\parallel_F,\label{theoremF}
\end{equation}
\begin{equation}
\mid\alpha_j-\beta_j\mid\leq\sigma_{\max}\left(\B-\A\right),\ j=1,2,\dotsc,n'.\label{abs}
\end{equation}
\end{thm}

\begin{pr} \upshape{First,} we take the T-SVD of $\A$ and $\B$ to get the f-diagonal tensors $\s_{\A}$ and $\s_{\B}$ with frontal slices
\begin{equation*}
\left(\sjj_{\A}\right)_i=\diag\left(\sigma^{i}_1,\dotsc,\sigma^{i}_{n'}\right),\ \left(\sjj_{\B}\right)_i=\diag\left(\tau^{i}_1,\dotsc,\tau^{i}_{n'}\right),
\end{equation*}
where the diagonal elements are arranged in descending order.
The Frobenius norm inequality of Mirsky type for matrices reads as
\begin{equation*}
\sum_{j=1}^{n'}\left(\sigma^i_{j}-\tau^i_j\right)^2\leq\parallel\Bjj_i-\Ajj_i\parallel_F^2.
\label{matrixF}
\end{equation*}
From (\ref{singualr}),
\begin{comment}
\begin{equation*}
\alpha_j=\frac{1}{r}\sum_{i=1}^r\sigma_j^i,\
%\alpha_1\geq\dotsc\geq\alpha_{n'}\geq 0,\
\beta_j=\frac{1}{r}\sum_{i=1}^r\tau_j^i\
%\beta_1\geq\dotsc\geq\beta_{n'}\geq 0,
\label{singularvalueAB}
\end{equation*}
\end{comment}
$\alpha_j=\frac{1}{r}\sum_{i=1}^r\sigma_j^i$ and $\beta_j=\frac{1}{r}\sum_{i=1}^r\tau_j^i$ are the $j$th singular values of $\A$ and $\B$ respectively.
Then, it is not hard to show
\begin{equation}
\begin{split}
\sum_{j=1}^{n'}\left(\alpha_j-\beta_j\right)^2=
\sum_{j=1}^{n'}\frac{1}{r^2}\left[\sum_{i=1}^r\left(\sigma_j^i-\tau_j^i\right)\right]^2
\leq\frac{1}{r}\sum_{j=1}^{n'}\sum_{i=1}^r\left(\sigma_j^i-\tau_j^i\right)^2\\
\end{split}.
\label{first}
\end{equation}
Adjusting the order of summation and regarding $\B-\A$ as one new tensor, we can get
\begin{equation}
\frac{1}{r}\sum_{i=1}^r\sum_{j=1}^{n'}\left(\sigma_j^i-\tau_j^i\right)^2\leq
\frac{1}{r}\sum_{i=1}^r\parallel\Bjj_i-\Ajj_i\parallel_F^2=
\frac{1}{r}\parallel\Bj-\Aj\parallel_F^2=\parallel\B-\A\parallel_F^2
\label{last},
\end{equation}
in which last equality follows from (\ref{Fnorm}). Finally, (\ref{theoremF}) is obtained by combining (\ref{first}) and (\ref{last}).

For (\ref{abs}), it is direct to show that
\begin{equation*}
\mid\alpha_j-\beta_j\mid=\mid\frac{1}{r}\sum^{r}_{i=1}\left(\sigma_j^i-\tau_j^i\right)\mid
\leq\frac{1}{r}\sum^{r}_{i=1}\mid\sigma_j^i-\tau_j^i\mid
\leq\frac{1}{r}\sum^{r}_{i=1}\sigma_{\max}\left(\Bjj_i-\Ajj_i\right),
\end{equation*}
where the second inequality is also from (\ref{matrixnorm}) and the right hand side of the equation is exactly equal to $\sigma_{\max}(\B-\A)$. That is, $\mid\alpha_j-\beta_j\mid\leq\sigma_{\max}\left(\B-\A\right),\ j=1,2,\dotsc,n'$.\qedsymbol
\end{pr}

%%===============第五部分：关于奇异值的：和或积有关性质======================
%%------顺序：先把加的2个和积的2个写出来：放在一个定理里-------------------
%%-------然后把涉及到正交张量的那个性质拿出来写一个定理---------------------
%%-------最后一个定理写扰动分析，最大最小边界------------------------------
\section{Tensor singular values of products and sums}
\hskip 2em The Weyl's theorem was extended from matrices to the T-product based tensors by Liu and Jin in \cite[Theorem 4.2]{Liu21}. 
Now, we are interested in developing useful inequalities for the singular values of the products and sums of tensors.

%%%===================两个和，两个积的===============================================
%%%================定理5.1：共计四个结论：只证前两个，后两个同理=====================
\begin{thm} Let $\A,\B\in\RR^{m \times n \times r}$ and $p=\min\{m,n\}$. The following inequalities hold for the decreasingly ordered singular values of $\A,\B,\A+\B$ and $\A*\B^{\T}$:
\begin{equation}%%%和
\sigma_{i+j-1}\left(\A+\B\right)\leq\sigma_i\left(\A\right)+\sigma_j\left(\B\right),
\label{he}
\end{equation}
\begin{equation}%%%积1
\sigma_{i+j-1}\left(\A*\B^{\T}\right)\leq r\sigma_i\left(\A\right)\sigma_j\left(\B\right),
\label{jione}
\end{equation}
\begin{equation}%%%差
|\sigma_{i}\left(\A+\B\right)-\sigma_i\left(\A\right)|\leq\sigma_1\left(\B\right).
\label{cha}
\end{equation}
\begin{equation}%%%积2
\sigma_{i}\left(\A*\B^{\T}\right)\leq r\sigma_i\left(\A\right)\sigma_1\left(\B\right),
\label{jitwo}
\end{equation}
For (\ref{he}) and (\ref{jione}), $1\leq i,j\leq p$ and $i+j-1\leq p$. For (\ref{cha}) and (\ref{jitwo}), $i=1,\dotsc,p$.
\end{thm}

\begin{pr}\upshape{As we can see from} \cite[Theorem 3.3.16]{Horn05}, if $r=1$ in our Theorem, all the conclusions are right. Thus, our proof is based on the theorem for matrix cases. Combining the FFT, T-SVD and Lemma 2.1, we can derive that
\begin{align*}
&\sigma_{i+j-1}\left(\A+\B\right)
=\frac{1}{r}\sum^{r}_{k=1}\sigma_{i+j-1}\left(\Ajj_k+\Bjj_k\right)\\
\leq&\frac{1}{r}\sum^{r}_{k=1}\sigma_{i}\left(\Ajj_k\right)+\sigma_{j}\left(\Bjj_k\right)
=\sigma_i\left(\A\right)+\sigma_j\left(\B\right),
\end{align*}
\begin{equation}
\sigma_{i+j-1}\left(\A*\B^{\T}\right)
=\frac{1}{r}\sum^{r}_{k=1}\sigma_{i+j-1}\left(\Ajj_k\Bjj_k^{*}\right)
\leq\frac{1}{r}\sum^{r}_{k=1}\sigma_{i}\left(\Ajj_k\right)\sigma_{j}\left(\Bjj_k\right).
\label{yiban}
\end{equation}
The proof of (\ref{jione}) is complete by exploiting below the Cauchy-Schwarz inequality for the right-hand side of (\ref{yiban}),
\begin{align*}
\sum^{r}_{k=1}\sigma_{i}\left(\Ajj_k\right)\sigma_{j}\left(\Bjj_k\right)
&\leq\sqrt{\sum^{r}_{k=1}\sigma_{i}\left(\Ajj_k\right)^2\sum^{r}_{k=1}\sigma_{j}\left(\Bjj_k\right)^2}\\
&\leq\sqrt{\left(\sum^{r}_{k=1}\sigma_{i}\left(\Ajj_k\right)\right)^2\left(\sum^{r}_{k=1}\sigma_{j}\left(\Bjj_k\right)\right)^2}\\
&=\sum^{r}_{k=1}\sigma_{i}\left(\Ajj_k\right)\sum^{r}_{k=1}\sigma_{j}\left(\Bjj_k\right)\\
&=r^2\sigma_i(\A)\sigma_j(\B).
\end{align*}
\begin{comment}
\begin{align*}
\sum^{r}_{k=1}\sigma_{i}\left(\Ajj_k\right)\sigma_{j}\left(\Bjj_k\right)
&\leq\sqrt{\sum^{r}_{k=1}\sigma_{i}\left(\Ajj_k\right)^2\sum^{r}_{k=1}\sigma_{j}\left(\Bjj_k\right)^2}\\
&\leq\sqrt{\left(\sum^{r}_{k=1}\sigma_{i}\left(\Ajj_k\right)\right)^2\left(\sum^{r}_{k=1}\sigma_{j}\left(\Bjj_k\right)\right)^2}\\
&=\sum^{r}_{k=1}\sigma_{i}\left(\Ajj_k\right)\sum^{r}_{k=1}\sigma_{j}\left(\Bjj_k\right)\\
&=r^2\sigma_i(\A)\sigma_j(\B).
\end{align*}
\end{comment}
The proofs of (\ref{cha}) and (\ref{jitwo}) are similar.\qedsymbol
\end{pr}

\begin{remark} \upshape{If} we let $i=j=1$ in (\ref{he}) and let $i=p$ in (\ref{cha}) respectively, it is clear that
\begin{equation*}
\sigma_{\max}(\A+\B)\leq\sigma_{\max}(\A)+\sigma_{\max}(\B),\ \sigma_{\min}(\A+\B)\geq\sigma_{\min}(\A)-\sigma_{\max}(\B).
\end{equation*}
\end{remark}

%%===================定理5.2：正交张量的那个=======================
\begin{thm} Given $\A\in\RR^{m \times n \times r}$ and two partially orthogonal tensors $\U\in\RR^{m \times k \times r}$ and $\V\in\RR^{n \times k \times r}$, where $k\leq\min\{m,n\}$. Then,
\begin{equation*}
\sigma_i\left(\U^{\T}*\A*\V\right)\leq\sigma_i(\A),\ i=1,\dotsc,k.
\end{equation*}
\end{thm}

\begin{pr} \upshape{If} $r=1$, this theorem is just the matrix case in \cite[Theorem 3.3.1]{Horn05}. It is easy to verify this inequality by the same technique in proving (\ref{he}).\qedsymbol
\end{pr}

%%=====================最后给一个扰动结果：双侧不等式=============================
\hskip 2em For a multiplicative perturbation to a tensor, we have the result below, which reduces to the matrix case for $r=1$.
\begin{thm} Let $\B'\!=\!\B\!+\!\delta\B\!=\!\U*\B*\V\!\in\!\RR^{m \times n \times r}$. Then,
\begin{equation*}
\sigma_{i}\left(\B'\right)
\leq r^2\sigma_{i}\left(\B\right)\sigma_{\max}\left(\U\right)\sigma_{\max}\left(\V\right).
\end{equation*}
\end{thm}

\begin{pr} \upshape{Since} $\B^{\T}$ has the same singular values with $\B$, then using (\ref{jitwo}), we get
\begin{align*}
\sigma_i(\B')=\sigma_i\left(\U*\B*\V\right)&=\sigma_i\left(\V^{\T}*\B^{\T}*\U^{\T}\right)
\leq r\sigma_i\left(\V^{\T}*\B^{\T}\right)\sigma_1(\U)\\
&=r\sigma_i(\B*\V)\sigma_1(\U)\leq r^2\sigma_i(\B)\sigma_1\left(\V^{\T}\right)\sigma_1(\U)\\
&=r^2\sigma_i(\B)\sigma_1\left(\V\right)\sigma_1(\U).\qedsymbol
\end{align*}
\end{pr}

\begin{comment}
\begin{remark} \upshape{Our Theorem is consistent with that in matrix case when} $r=1$.
\end{remark}
\end{comment}

%%%==========================最后一部分：总结===========================
\section{Concluding remarks}
\hskip 2em In this paper, we present some perturbation results for the T-product based tensor singular values, 
including the relationship between the singular values of the subtensor and the original tensor, the sums and products of tensors. 
The classical results given by Mirsky are also extended to tensors. Extensions of some other relative perturbation results for matrix singular values need to be further studied.

%%%=======================参考文献==================================


\begin{thebibliography}{99}
\bibliographystyle{plain}
\vspace{-1.0em}

%%%--------------------【1】：特征值，特征向量的论文=========================
\bibitem{Braman10}
        K. Braman,
        \emph{Third-order tensors as linear operators on a space of matrices},
        Linear Algebra Appl., 433 (7) (2010) 1241--1253.

%%%========================【2】矩阵计算第四版英文书======================================
\bibitem{Golub13}
        G. H. Golub, C. F. Van Loan,
        \emph{Matrix Computations},
        Johns Hopkins Univ. Press, $4$th edition, 2013.


%%%===================【3】Topic书===============================
\bibitem{Horn05}
        R. A. Horn, C. R. Johnson,
        \emph{Topics in Matrix Analysis},
        Cambridge University, 1991.

%%%=====================【4】2008年首次提出新型T-积定义===================
\bibitem{Kilmer08}
        M. E. Kilmer, C. D. Martin, L. Perrone,
        \emph{A third-order generalization of the matrix SVD as a product of third-order tensors},
        Tech. Report TR-2008-4, Tufts University, Computer Science Department, 2008.

%%%=====================【5】应用
\bibitem{Kilmer13}
         M. E. Kilmer, K. Braman, N. Hao, R. Hoover,
         \emph{Third-order tensors as operators on matrices: A theoretical and computational framework with applications in imaging},
         SIAM J. Matrix Anal. Appl., 34 (2013) 148--172.

%%%======================【6】新的奇异值定义：时间在后=================
\bibitem{Ling21}
        C. Ling, J. Liu, C. Ouyang, L. Qi,
        \emph{ST-SVD factorization and s-Diagonal tensors},
        arxiv:2104.05329, 2021.

%%%========================【7】张量特征值，不等式====================
\bibitem{Liu21}
         W. Liu, X. Jin,
        \emph{A study on T-eigenvalues of third-order tensors},
        Linear Algebra Appl., 612 (2021) 357--374.

%%%===============【8】2020年IEEE论文：本文定义============
\bibitem{Lu20}
        C. Lu, J. Feng, Y. Chen, W. Liu, Z. Lin, S. Yan,
        \emph{Tensor robust principal component analysis with a new tensor nuclear norm},
         IEEE Trans. Pattern Anal. Mach. Intell., 42 (2) (2020) 925--938.

%%%======================【9】应用
\bibitem{MiaoQi20}
         Y. Miao, L. Qi, Y. Wei,
         \emph{Generalized tensor function via the tensor singular value decomposition based on the T-product},
         Linear Algebra Appl., 590 (2020) 258--303.

%%%========================【10】Mirsky==============================
\bibitem{Mirsky60}
         L. Mirsky,
         \emph{Symmetric gauge functions and unitarily invariant norms},
         Quart. J. Math., 11 (1960) 50--59.

%%%=======================【11】新的奇异值定义：时间在前===============
\bibitem{Qi21}
         L. Qi, G. Yu,
         \emph{T-singular values and T-Sketching for third order tensors},
         arXiv:2013.00976, 2021.

%%%=====================【12】应用
\bibitem{Semerci14}
         O. Semerci, N. Hao, M. E. Kilmer, E. L. Miller,
         \emph{Tensor-based formulation and nuclear norm regularization for multienergy computed tomography},
         IEEE Trans. On Image Processing, 23 (2014) 1678--1693.

%%%=====================【13】1972矩阵交错定理============================
\bibitem{Thompson71}
         R. C. Thompson,
         \emph{Principal Submatrices IX: Interlacing Inequalities for Singular Values of Submatrices},
         Linear Algebra Appl., 5 (1) (1971) 1--12.

%%%=====================【14】应用
\bibitem{XiaoChen21}
        X. Xiao, Y. Chen, Y. J. Gong, Y. Zhou,
        \emph{Low-rank reserving t-linear projection for robust image feature extraction}, IEEE Trans. On Image Processing, 30 (2021) 108--120.

%%%=====================【15】应用
\bibitem{Yang16}
        L. Yang, Z. H. Huang, S. Hu, J. Han,
        \emph{An iterative algorithm for third-order tensor multi-rank minimization}, Comput Optim Appl., 63 (2016) 169--202.

%%%=====================【16】应用
\bibitem{Zhang17}
         Z. Zhang, S. Aeron,
         \emph{Exact tensor completion using t-SVD},
         IEEE Trans. on Signal Processing, 65 (2017) 1511--1526.

%%%=====================【17】应用
\bibitem{Zhou18}
        P. Zhou, C. Lu, Z. Lin, C. Zhang,
        \emph{Tensor factorization for low-rank tensor completion},
        IEEE Trans. on Image Processing, 27 (2018) 1152--1163.

\vspace{-1.0em}
\end{thebibliography}
\end{document}